\documentclass{article}
\usepackage[hmargin = 3.7 cm,vmargin = 3.7 cm]{geometry}
\begin{document}

\title{Volume functionals on pseudoconvex hypersurfaces}
\author{Author}
\date{\today}
\author{Simon Donaldson and Fabian Lehmann}
%\date{\today}
\maketitle
\newcommand{\bC}{{\bf C}}
\newcommand{\bR}{{\bf R}}
\newcommand{\bP}{{\bf P}}
\newcommand{\oPsi}{\overline{\Psi}}
\newcommand{\oz}{\overline{z}}
\newcommand{\db}{\overline{\partial}}
\newcommand{\cF}{{\cal F}}
\newcommand{\cU}{{\cal U}}
\newcommand{\cL}{{\cal L}}
\newcommand{\cH}{{\cal H}}
\newcommand{\cM}{{\cal M}}
\newcommand{\bE}{{\bf E}}
\newtheorem{lem}{Lemma}
\newtheorem{prop}{Proposition}

\

\

\     {\it Dedicated to Oscar Garc{\'i}a-Prada,  in celebration of his sixtieth birthday.}

\

\

\

The focus of this paper is on a volume form defined on a pseudoconvex hypersurface $M$ in a complex Calabi-Yau manifold (that is, a complex  $n$-manifold  with a nowhere-vanishing holomorphic $n$-form). We begin by  defining this volume form and observing that it can  be viewed as a generalisation of the affine-invariant volume form on a convex hypersurface in $\bR^{n}$. We compute the first variation, which leads to a similar generalisation of the affine mean curvature. In the second section we investigate the constrained variational problem, for pseudoconvex hypersurfaces $M$ bounding  compact domains $\Omega\subset Z$. That is, we  study critical points of the volume functional $A(M)$ where the ordinary volume $V(\Omega)$ is fixed. The critical points are analogous to constant mean curvature submanifolds.  We find that Sasaki-Einstein hypersurfaces satisfy the condition, and  in particular the standard sphere $S^{2n-1}\subset \bC^{n}$ does. The main work in the paper comes in the third section where  we compute the second variation about the sphere. We find that it is negative in \lq\lq most'' directions but non-negative in directions corresponding to
deformations of $S^{2n-1}$ by holomorphic diffeomorphisms. We are lead to conjecture a  \lq\lq minimax''characterisation of the sphere.   We also discuss connections with the affine geometry case and with K\"ahler-Einstein geometry.

Our original motivation for investigating these matters came from the case $n=3$ and the embedding problem studied in our previous paper \cite{kn:DL}. There are some special features in this case. The volume functional  can be defined without reference to the embedding in $Z$ using only a closed \lq\lq pseudoconvex'' real $3$-form  on $M$. We develop some of the theory from the point of the symplectic structure on exact $3$-forms on $M$ and the moment map for the action of the diffeomorphisms of $M$.

\

This work was supported by the Simons Foundation through the {\it Simons Collaboration on special holonomy in geometry, analysis and physics}. The authors are grateful to Erwin Lutwak for references concerning the affine isoperimetric inequality. 

\section{ The hypersurface volume form}

Let $Z$ be an $n$-dimensional complex manifold,  $\mu\in \Omega^{2n}(Z)$  be a volume form on $Z$ and $M\subset Z$  a pseudoconvex hypersurface. (Note that in this paper we shall mean by \lq\lq pseudoconvex''  what is sometimes called \lq\lq strongly pseudoconvex''). We define a volume form on $M$ as follows.  Let $F$ be a local defining function for $M=F^{-1}(0)$ near $p$, with orientation conventions chosen so that in the model case of $S^{2n-1}$ the function $F$ is positive on the exterior of the ball. Then we have a positive $2n$-form
  $ i\partial F\wedge \db F \wedge (i\partial\db F)^{n-1} $ which we can write as $ 2^{n} (n-1)! g \mu$ for a function $g$. At a point $p\in M$ let $w\in TZ_{p}$  be a vector such that $TZ_{p}= \bR w \oplus TM_{p}$. Now we have a volume form on $TM_{p}$ defined by the expression
$$  \nu_{p}=  2 g^{1/n+1}\ (\nabla_{w}F)^{-1}\ i_{w}(\mu). $$
(The factors are chosen so that this gives the standard volume form on $S^{2n-1}$, when $\mu$ is the standard volume form on $\bC^{n}$.)
\begin{lem}
$\nu_{p}$ is independent of the choice of $w$ and $F$
\end{lem}

Changing $w$ by the addition of a tangent vector to $M$ does not change
$\nabla_{w} F$ or $i_{w}\mu\vert_{M}$ and changing $w$ to $c w$ multiplies both terms by $c$. So the expression is independent of $w$. Changing $F$ to $cF$ changes $g$ to $c^{n+1} g$, so does not change $\nu_{p}$. If $h$ is a smooth function on $Z$ equal to $1$ on $M$ then
$$ \partial \db(hF)= (\partial \db h) F + \partial h \db F  - \db h \partial F + h \partial \db F  $$
and  at the point $p$ we have $\partial(hF)= h \partial F= \partial F$. So, at this point
$p$,  we have $ \partial (hF)\wedge \db (hF) \wedge (\partial\db (h F))^{n-1}= \partial F\wedge \db F \wedge (\partial\db F)^{n-1}$   and there is no change in $g$.

Another way of expressing the definition is that the {\it Levi form} $L$ of the hypersurface is invariantly defined in $\Lambda^{2} H^{*} \otimes (TM/H)$ where $H\subset TM$ is the complex tangent space $TM\cap ITM$. Then $L^{n-1}$ is defined in $\Lambda^{2n-2} H^{*} \wedge (TM/H)^{n-1}$. On the other hand $\mu$ is an element of $\Lambda^{2n-2} H^{*}\otimes \Lambda^{2}(TZ/H)^{*}$. The complex structure gives an isomorphism $TZ/TM= TM/H$ and so  $\Lambda^{2}(TZ/H)^{*}= (TM/H)^{-2}$. Thus  the ratio $L^{n-1}/ \mu$ can be regarded as an  element of the line
$(TZ/TM)^{n+1}$. Then
   $$  \left( L^{n-1}/\mu\right)^{1/n+1}\in TZ/TM $$  
   and the contraction of $\mu$ with this element defines $\nu\in \Lambda^{2n-1}T^{*}M$.
   
   Another variant of the definition  will be useful. Given a defining function $F$ choose
a normal vector $w$ such that $\nabla_{w} F=1$. Then we can produce two 
 volume forms on $M$
\begin{equation}   \nu_{1}= i_{w} \mu \ \ \ \ , \ \ \ \ \nu_{2} =\frac{1}{(n-1)!} (i\partial \db F)^{n-1} I d F
. \end{equation}
If we multiply $F$ by $c$ then $\nu_{1}$ multiplies by $c^{-1}$ and $\nu_{2}$
by $c^{n}$. So the combination $$
\nu_{2}^{1/n+1} \nu_{1}^{n/n+1} $$ is independent of the choice of $F$ and
this is our volume form $\nu$. (Note that if $\sigma_{1}, \sigma_{2}$ are volume forms on a manifold and $a+b=1$ then $\sigma_{1}^{a}\sigma_{2}^{b}$ is a volume form.)

   To make all this more explicit, suppose that $\mu$ is the standard volume form in local complex co-ordinates $z_{1}, \dots ,z_{n}$ (which is our main case of interest). Writing $z_{a}= x_{a}+ i y_{a}$, suppose that $M$ is the graph of a function $x_{n}= f(z_{1}, \dots, z_{n-1}, y_{n})$, so we can regard
$z_{1}, \dots, z_{n-1}, y_{n}$ as local coordinates on $M$. Take
$F= x_{n}- f$ and $w$ to be the constant vector $\frac{\partial}{\partial x_{n}}$. So the volume form $\nu_{1}$ is the standard form in these coordinates.
In the calculation of $\nu_{2}$ we have 
$ (i\partial \db F)^{n-1}= \sum G^{ab} \sigma_{ab}$ where
$\sigma_{ab}$ is the $2n-2$ form given by contracting the volume form of $\bC^{n}$ with $\frac{\partial}{\partial z_{a}}$ and $\frac{\partial}{\partial \overline{z}_{b}}$ and $(G^{ab})$ is the matrix of cofactors of  $ \frac{\partial^{2} f}{\partial
z_{a} \partial\overline{z}_{b}}$. 
To compute $\nu_{2}$ we take the wedge product of $\sum G^{ab}\sigma_{ab}$ with the $1$-form $$IdF = -dy_{n} + \frac{\partial f}{\partial y_{n}} dx_{n} +\sum_{a=1}^{n-1}
 \left( \frac{\partial f}{\partial y_{a}}dx_{a} - \frac{\partial f}{\partial x_{a}} dy_{a}\right),$$
and finally substitute $$dx_{n}= \frac{\partial f}{\partial y_{n}} dy_{n} + \sum_{a=1}^{n-1} \left(\frac{\partial f}{\partial x_{a}} dx_{a} + \frac{\partial f}{\partial y_{a}} dy_{a}\right). $$

The calculation is much clearer at a point where the derivative of $f$ vanishes. At such a point we have the  simple formula
$ \nu= J^{1/n+1} dx_{1}dy_{1}  \dots dy_{n}$
where $J$ is the determinant of the $(n-1)\times (n-1)$ matrix  $ \left(\frac{\partial^{2}
f}{\partial
z_{a} \partial\overline{z}_{b}}\right)_{a,b=1\dots, n-1}$.

\

The volume form construction can be viewed as a generalisation of the \lq\lq affine volume'' $\nu_{{\rm aff}}$ on a convex hypersurface $S$ in $\bR^{n}$, invariant under $SL(n,\bR)$. The second fundamental form $B$  of $S$ is defined affine-invariantly in $s^{2}(T^{*}S)\otimes N$ where $N$ is the normal bundle $\bR^{n}/TS$. Then the determinant of $B$ is defined in the line $(\Lambda^{n-1}T^{*}S)^{2}\otimes N^{n-1}$. We have an isomorphism $\Lambda^{n-1} T^{*}S \otimes N^{*}= \left(\Lambda^{n}\bR^{n}\right)^{-1}$
and since $\Lambda^{n}\bR^{n}$ has a fixed trivialisation we get $\Lambda^{n-1}T^{*}S= N$. So ${\rm det}\ B$ lies in $ (\Lambda^{n-1}T^{*} S)^{n+1}$ and $( {\rm det}\ B)^{1/n+1}$ is a volume form on $S$. Now regard $\bR^{n}$ as embedded in $\bC^{n}$ in the standard way and define $M= S\times i \bR^{n}\subset \bC^{n}$. This is a pseudoconvex hypersurface and tracing through the definitions we find that the constuctions match up: {\it i.e.} the volume form $\nu$ on $M$ is the product of the affine volume on $S$ and the standard volume on $i \bR^{n}$. 

\

To see this more explicitly, let $S$ be described locally as a graph
$ x_{n}= f(x_{1}, \dots, x_{n-1})$ of a convex function $f$ on $\bR^{n-1}$. Then $M$ is described by the same function $f$, independent of the imaginary coordinates $y_{a}$. The matrix  $ \frac{\partial^{2} F}{\partial z_{a} \partial\overline{z}_{b}}$. has the block form
$$  \left( \begin{array}{cc}   \frac{\partial^{2} f}{\partial x_{a} \partial x_{b}} &0\\0&0 \end{array} \right). $$
There is only one non-zero cofactor and we again get the simple formula
$$  \nu= \left( {\rm det}(\frac{\partial^{2} f}{\partial x_{a} \partial
x_{b}} )\right)^{1/n+1}  dx_{1}\dots dx_{n-1} dy_{1}\dots dy_{n},  $$
which agrees with standard formula for the affine volume of a graph.

\section{The first variation}

From now on we consider the volume form $\mu$ on $Z$ arising from a complex Calabi-Yau structure----a non-vanishing holomorphic $n$-form $\Psi$. Thus $$   \mu= c_{n} \Psi \wedge\oPsi, $$
where $c_{n}=  2^{-n} (i)^{n^{2}}$ is the constant such that
$   c_{n} dz_{1}\dots dz_{n} d\oz_{1}\dots d\oz_{n}$ is the standard volume form on $\bC^{n}$.
Suppose  that $M\subset Z$ is a compact pseudoconvex hypersurface, so we have a total volume
  $$   A(M)= \int_{M} \nu. $$
  We want to compute the first variation of this functional. To set this up we review some structure theory, as in \cite{kn:DL}.  The restriction of $\Psi$ to $M$ can be written as $\theta \wedge \chi$ where $\chi=\alpha+ i \beta$ for real $(n-1)$ forms $\alpha,\beta$ and $\theta$ is a contact $1$-form for the contact structure $H$.
If $\tilde{\theta}=I \theta$ then at a point of $M$ the holomorphic $n$-form is $(\theta+i\tilde{\theta})\wedge \chi$.  Just as in \cite{kn:DL} we can fix the choice of $\theta$ by the condition that 
 $$    \theta\wedge (d\theta)^{n-1}= (n-1)!\ c_{n-1} \theta \wedge \chi\wedge \overline{\chi}. $$
Then the volume form is
$$  \nu= c_{n-1} \theta\wedge \chi\wedge\overline{\chi}, $$
(which is also $ 1/(n-1)! \theta\wedge (d\theta)^{n-1}$, by our normalisation.
 When $n$ is odd we can write $\nu$ in terms of either $\alpha$ or $\beta$:
$$ \nu = 2 c_{n-1} \theta \wedge \alpha^{2}= 2 c_{n-1}\theta\wedge\beta^{2}$$ and $\alpha\wedge\beta=0$.  When $n$ is even
$$  \nu = 2 i c_{n-1}\theta \wedge \alpha\wedge \beta. $$

 The choice of $\theta$ defines a Reeb vector field $v$ transverse to $H$ with $i_{v}(d\theta)=0$ and $\theta(v)=1$. The choice of $\chi$ is fixed by requiring that $i_{v}(\chi)=0$. 

More generally if $\Psi\vert_{M}= \theta\wedge \chi$ and we do not normalise the choice of contact form $\theta$ as above then we have a formula
\begin{equation}  \nu= \left( c_{n-1} \theta\wedge \chi\wedge \overline{\chi}\right)^{n/n+1}\ \left( \frac{1}{(n-1)!} \theta\wedge (d\theta)^{n-1}\right)^{1/n+1}. \end{equation}

(So the normalisation is just the requirement that

\

Now suppose that we have a $1$-parameter deformation $ m_{t}:M\rightarrow Z$ with $m_{0}$ the given embedding. The derivative at $t=0$ is a section $w$ of $TZ\vert_{M}$. 

\begin{prop}

The derivative of $A$ is, for $n$ odd,
\begin{equation} \delta A = -2^{1-n}\frac{n}{n+1}  {\rm Re} \left( \int_{M} (i_{w}\Psi) \wedge d\overline{\chi}\right), \end{equation} and for $n$ even
$$ \delta A= -2^{1-n} \frac{n}{n+1} {\rm Im}  \left( \int_{M} (i_{w}\Psi)
\wedge d\overline{\chi}\right). $$
\end{prop}

 For simplicity we will write the proof for the case $n$ odd. 
The derivative of $m_{t}^{*}(\Psi)$ is given by the Cartan formula
$$  \delta \Psi = d(i_{w}\Psi) + i_{w} d\Psi= d (i_{w} \Psi) . $$
 Consider first the case when $m_{t}$ are diffeomorphisms to the same image $M$, so $w$ is tangent to $M$. Clearly $\delta A=0$, so we want to check that the right hand side of (3) also vanishes.

Restricting to $M$, we have  $\Psi= \theta \wedge \chi$ so
$$ i_{w}(\Psi) = i_{w}\theta \wedge \chi -  \theta \wedge i_{w}(\chi). $$
The fact that $d\Psi\vert_{M}=0$ implies that
$$  d{\chi}= L_{v} \chi \wedge \theta$$
where $L_{v}$ is the Lie derivative along the Reeb vector field $v$. So
$$   (i_{w}\Psi) \wedge d\overline{\chi}= \theta(w) \chi\wedge d\overline{\chi}, $$
and the real part is $\theta(w) d(\chi\wedge \overline{\chi})$. But 
$$  d(\chi\wedge \overline{\chi})= L_{v}(\chi\wedge \overline{\chi})\wedge \theta = L_{v}(d\theta^{n-1})\wedge \theta=0, $$
confirming the desired vanishing of the right hand side of (3), in this case.

Next consider the case when $m_{t}$ preserves the contact structure $H$, so  $\delta \Psi$ has the form $\theta \wedge \lambda$. The formula (2) for the volume form shows that
$$  \delta \nu= 2^{1-n}\frac{n}{n+1} {\rm Re}\ \left( \lambda\wedge \overline{\chi} \wedge \theta\right), $$
which is equal to the real part of $$  2^{1-n} \frac{n}{n+1}\delta\Psi \wedge  \overline{\chi}.$$

Now we have $\delta \Psi= d(i_{w} \Psi)$ so
$$  \delta A= \int_{M} \delta \nu=  2^{1-n}\frac{n}{n+1} {\rm Re} \int_{M} d(i_{w} \Psi) \wedge  \overline{\chi}. $$
Applying Stokes' Theorem, this is:
$$\delta A= - 2^{1-n} \frac{n}{n+1} {\rm Re} \int_{M} i_{w} \Psi \wedge d\overline{\chi}, $$
as required

Any infinitesimal variation of a contact structure can be realised by the action of a vector field. (This is
the infinitesimal form of Gray's theorem on the stability of contact structures.) It follows from this that the formula (3)  holds for general variations.

\

Note that the fact that $\chi\wedge\theta$ is closed implies that  $d\overline{\chi} = L_{v}\overline{\chi} \wedge \theta$.

\

For our hypersurface $M\subset Z$ we have a distinguished normal vector field $Iv$, so it suffices to consider variations with $w= f (Iv)$, for a function $f$ on $M$. In that case $i_{w}\Psi= i f \chi$ and we get
$$  \delta A= \int_{M} f \ {\bf h}\  \nu, $$
where, in the case $n$ odd,
$$  {\bf h}\  \nu =  2^{1-n}\frac{n}{n+1}\  {\rm Im}( \chi \wedge L_{v}\overline{\chi})\wedge \theta. $$ In terms of the real and imaginary parts $\chi=\alpha+i \beta$ this is
$$  {\bf h}\ \nu=  -2^{2-n} \frac{n}{n+1} \beta\wedge L_{v}\alpha \wedge \theta. $$
For $n$ even the formula is
$$  {\bf h}\ \nu= -2^{1-n} \frac{n}{n+1} (\alpha\wedge L_{v}\alpha +\beta\wedge L_{v}\beta) \wedge \theta. $$

This function ${\bf h}$ is an invariant of the pseudoconvex hypersurface in the complex Calabi-Yau manifold $(Z,\Psi)$ and is the analogue of the mean curvature in our setting. In the case when $M=S\times i\bR^{n}$ it reduces to  the classical  affine mean curvature of a convex hypersurface.
There are also analogues,  {\bf h}={\rm const}, of constant mean curvature hypersurfaces.  In the case when $M$ is the boundary of a compact pseudoconvex domain $\Omega \subset Z$ we have another functional 
$$   V(M) = {\rm Vol}(\Omega)= \int_{\Omega} c_{n} \Psi \wedge \overline{\Psi}. $$
Under our variation field $f (Iv)$ the variation of $V$ is
simply
$$  \delta V= \int_{M} f \nu, $$
so the equation ${\bf h}= {\rm constant}$ characterises the extrema of $A(M)$ with $V(M)$ fixed.

 A submanifold $M$ satisfying
$$   L_{v} \chi = \lambda i\  \chi, $$
for constant $\lambda>0$ is called {\it Sasaki-Einstein}. Clearly a  Sasaki-Einstein submanifold has  
  ${\bf h}$ constant.  We have a partial converse for this, in the case $n=3$. In that case, for a compact pseudoconvex $M^{5}\subset Z$, we defined in \cite{kn:DL} a finite dimensional space ${\cal H}$ which arises as the obstruction space for the deformation problem studied in \cite{kn:DL}. This is the kernel of the operator
   $$d_{H}: \Omega^{-}_{H}\rightarrow \Omega^{3}_{H}.$$
   Here $\Omega^{-}_{H}$ denotes the space of anti-self-dual $2$-forms on $H\subset TM$---the orthogonal complement of $\alpha,\beta, d\theta$---and $d_{H}$ is the horizontal component of the exterior derivative. We can write ${\bf h}$ as $L_{v}\alpha\wedge \beta$ (which is also $-\alpha \wedge L_{v}\beta$, since $\alpha\wedge \beta=0$).

\begin{prop}
 Suppose that $M^{5}\subset Z$ has ${\cal H}=0$ and that ${\bf h}$ is a positive constant $\mu$. Then $M$ is Sasaki-Einstein.
\end{prop}
In the proof we will use various identities for which we refer the reader to \cite{kn:DL}.  Recall that 
$$  {\bf h}\ \nu= \frac{3}{8} L_{v} \alpha\wedge \beta\wedge \theta = -\frac{3}{8}  \alpha\wedge L_{v} \beta\wedge \theta. $$

 Let $\gamma_{1}= L_{v}\alpha- (8\mu/3) \beta$. Then 
$$  d_{H}(\gamma_{1})= d_{H}(L_{v}\alpha) - (8\mu/3) d_{H}\beta= L_{v}(d_{H}\alpha) -(8\mu/3) d_{H}\beta =0, $$
since $d_{H}\alpha=d_{H}\beta=0$. We claim that $\gamma_{1}$ is an anti-self-dual form. For the assumption ${\bf h}=\mu$ gives
$ L_{v}\alpha\wedge \beta=(8\mu/3) \beta^{2}$ so  $\gamma_{1}\wedge\beta$ vanishes. The fact that $ L_{v}\alpha\wedge \alpha=0$ gives $\gamma_{1}\wedge \alpha=0$ and similarly for $\gamma_{1}\wedge d\theta$. Thus $\gamma_{1}$ lies in ${\cal H}$. In the same way, one sees that $\gamma_{2}= L_{v}\beta+ (8\mu/3) \alpha$ lies in ${\cal H}$. Under the hypotheses of the proposition $\gamma_{1}, \gamma_{2}$ both vanish, which means that $M$ is Sasaki-Einstein (with $\lambda= 8\mu/3$).

\

The $5$-manifold  $S^{2}\times T^{3}$, thought of as a hypersurface in $\bC^{3}/ i {\bf Z}^{3}$, is an example which is not Sasaki-Einstein but with ${\bf h}$ a positive constant.  

\

\section{The second variation}

The scaling behaviour of the functional $A(M)$ for submanifolds $M\subset \bC^{n}$ is
$$  A(\kappa M)= \kappa^{2n^{2}/n+1} A(M), $$
so the quantity
$$  R(M)= \frac{A(M)}{\ \ \ \ \  V(M)^{n/n+1}}$$
is scale-invariant. The critical points of the functional $R$ are the manifolds with constant ${\bf h}$, considered above. The analogue in the affine-geometry case, for a convex hypersurface $S\subset \bR^{n}$ bounding a domain $\Omega_{\bR}$ is the ratio
$$   R_{{\rm aff}}(S)= \frac{A_{{\rm aff}}(S)}{\ \ \ \ \  \left( {\rm Vol}\  \Omega_{\bR}\right)^{n-1/n+1}}. $$
The \lq\lq affine isoperimetric inequality''  asserts that this functional $R_{{\rm aff}}$ is maximised by the ellipsoids. This was proved by Blaschke for $n=3$ in \cite{kn:Blaschke} and extended to higher dimensions by several authors \cite{kn:Deicke}, \cite{kn:Nakajima},\cite{kn:Santalo}. Thus  it is natural to ask if there could be some similar statement in the complex case, for the functional $R$. In that direction, we will now compute the second variation of $R$ at the standard sphere $S^{2n-1}\subset \bC^{n}$.

For the calculation we use our third description of the volume form 
$\nu = \nu_{2}^{1/n+1} \nu_{1}^{n/n+1}$. We use the standard decomposition $TS^{2n-2}= H \oplus  \bR v$ with the Reeb vector field $v$ and corresponding contact form $\theta$. We write $d_{H}$ for the $H$-component of the exterior derivative.  In polar coordinates the tangent bundle of $\bC^{n}\setminus \{0\}$ is represented (with a slight stretch of notation) as $H\oplus \bR v \oplus \bR \partial_{ r}$: the complex structure is given by
$$   I_{\bC^{n}} \partial_{r} = -r^{-1} v\ \ ,\ \  I_{\bC^{n}} v = r \partial_{r}  , $$
and the fixed complex $I$ structure on $H$.

  We consider a hypersurface $M_{f}$ defined  in polar co-ordinates by an equation $r= e^{f}$ where $f$ is a function on $S^{n-1}$. There is an obvious radial projection map $m_{f}: S^{2n-1}\rightarrow M_{f}$. A defining  function in the sense of Section 1 is $F=
r-e^{f}$. We take $w$ to be the vector field $\partial_{r}$ so $\nabla_{w}F=1$.
In Section 1 we defined, in this situation, volume forms $\nu_{1}, \nu_{2}$ on $M_{f}$. Write
$\nu_{1,f}= m_{f}^{*}(\nu_{1})\ ,\  \nu_{2,f}= m_{f}^{*}(\nu_{2})$. These are volume forms on $S^{2n-1}$. We have

\begin{equation}   A(M_{f})= \int_{S^{2n-1}}  \nu_{1,f}^{n/n+1} \nu_{2,f}^{1/n+1}. \end{equation}

The calculation of the volume form $\nu_{1,f}$ is straightforward. It is just
$$   \nu_{1, f}= e^{(2n-1) f}  {\rm vol}, $$
where ${\rm vol}$ is the standard volume form on the sphere. 
  The volume form $\nu_{2,f}$ is 
\begin{equation}  \nu_{2,f}= (2^{n-1}(n-1)!)^{-1}\ \ \gamma \wedge d\gamma^{n-1}\end{equation}
 where $\gamma= m_{f}^{*}(I_{\bC^{n}} dF)$.
 
 We have
$$  dF= dr- e^{f} df= dr - e^{f}(d_{H} f + L_{v} f \theta). $$
Now $I_{\bC^{n}} dr=-r\theta$ and $I_{\bC^{n}}\theta= dr/r$. So
$$  I_{\bC^{n}} dF = -r \theta - e^{f}( I d_{H} f + (L_{v} f) r^{-1} dr). $$ 
From this we get
\begin{equation} -\gamma= e^{f}\left( (1+(L_{v} f)^{2}) \theta + (L_{v} f) d_{H} f + I d_{H} f\right) . \end{equation}

Substituting (6) into (5) we get a somewhat unwieldy  formula for $\nu_{2,f}$
in terms of the derivatives of $f$ on the sphere. Combining with the simple formula  for $\nu_{1,f}$ we arrive at an explicit formula  for $A(M_{f})$. The volume $V(M_{f})$ is simply
$$   V(M_{f})= \frac{1}{2n} \int_{M} e^{2nf} \ {\rm vol}, $$
so we can write down a formula for  $R(M_{f})=A V^{-n/n+1}$.

For our purposes here we only need the second variation about $f=0$ so the formulae become more tractable. We replace $f$ by $tf$, with $f$ fixed  and a real parameter $t$. We know that $R_{t}=R(M_{tf})$ is constant to first order and we are concerned with the second order term
$$  R_{t}= R(S^{2n-1})(1 + t^{2} Q(f) + O(t^{3}), $$
where the expression $Q(f)$ is quadratic in $f$ and its derivatives. Our problem is to compute $Q$. We expand the expression (6) for $\gamma$ up to order $t^{2}$:
$$   -\gamma= \theta + t\left(f \theta + I d_{H} f\right) + t^{2}\left( (f^{2}/2) \theta + (L_{v}f)^{2}\theta + f I d_{H} f+ L_{v}f d_{H} f\right) + O(t^{3}), $$
and then get an expansion of $\gamma\wedge d\gamma^{n-1}$ up to order $t^{2}$ and so on. The expressions are still complicated, but one sees from an outline calculation (and  applications of Stokes' Theorem) that in the end the quadratic form $Q$ is some linear combination of five constituents:
$$I(f)^{2}\ \ ,\ \  \Vert f\Vert^{2}\ \ ,\ \  \Vert d_{H}f\Vert^{2}\ \ ,\ \  \Vert L_{v} f\Vert^{2}\ \ ,\ \  \Vert \Delta_{H} f\Vert^{2}. $$
Here $I(f)$ denotes the integral of $f$ with respect to the standard volume form ${\rm vol}$, divided by ${\rm Vol}\ (S^{2n-1})^{1/2}$. The norm $\Vert\ \ \Vert$ is the standard $L^{2}$ norm. The operator $\Delta_{H}$ is $d_{H}^{*}d_{H}$, so
$\langle \Delta_{H} f_{1}, f_{2}\rangle = \Vert d_{H} f\Vert^{2}$. 
(The $\Vert \Delta_{H}f\Vert^{2} $ contribution arises in the following way. We have $4 (n-2)!\Delta_{H} f vol= d_{H}I d_{H} f \wedge (d\theta)^{n-2}\wedge \theta$ and this appears in the $O(t)$ term of $\nu_{2}$. Then $(\Delta_{H} f)^{2}$ appears in the $O(t^{2})$ term of $\nu_{2}^{1/n+1}$.)

\

The fact that $Q$ has this general shape, combined with the invariance properties of the functionals, determine $Q$, up to an overall factor. If $f$ is a constant function then $M_{f}$ is another round sphere and the scale invariance means that $R(M_{f})= R(S^{2n-1})$. So $Q$ vanishes on the constant functions and we have
\begin{equation}  Q= a_{0}(\Vert f\Vert^{2}- I(f)^{2})+ a_{1}\Vert d_{H}f\Vert^{2}+ a_{2} \Vert L_{v} f\Vert^{2}+ a_{3} \Vert \Delta_{H} f\Vert^{2}, \end{equation}
for some constants $a_{0},a_{1},a_{2}, a_{3}$. 

We apply some harmonic analysis.
The  functions on $S^{2n-1}$ decompose into eigenspaces for the Laplacian, which correspond to harmonic polynomials on $\bR^{2n}= \bC^{n}$. A harmonic polynomial of degree $d$ gives an eigenfunction with eigenvalue $d(d+2n-2)$. It is simpler here to work with complex-valued functions. Let $s^{p,q}$ be the space of complex polynomials $P(z,\overline{z})$ on $\bC^{n}$ of degree $p$ in $z$ and $q$ in $\overline{z}$. Let $h^{p,q}\subset s^{p,q}$ be the kernel of the contraction map $s^{p,q}\rightarrow s^{p-1, q-1}$. The harmonic polynomials of degree $d$ are the sum of the $h^{p,q}$ with $p+q=d$. The standard circle action (multiplication by unit complex numbers) on the sphere acts with   with weight $p-q$ on $h^{p,q}$. 

The spaces $h^{p,1}$ have special significance in our situation. A element of $s^{p,1}$ can be identified, using the metric on $\bC^{n}$, with a holomorphic vector field: we map
$ P=\sum P_{a}(z) \overline{z}_{a}$ to the vector field
$$ V= \sum P_{a}(z) \frac{\partial}{\partial z_{a}}. $$
The contraction $s^{p,1}\rightarrow s^{p-1,0}$ gives the divergence
$$  \sum \frac{\partial P_{a}}{\partial z_{a}}. $$
So the $h^{p,1}$ correspond to holomorphic volume-preserving vector fields with polynomial co-efficients. Let $\xi$ be  point on the unit sphere. The normal component of the above vector field $V$ at $\xi$ is
$${\rm Re}( \sum P_{a}(\xi) \overline{\xi}_{a}),$$ 
which is just the real part of the value at $\xi$ of the  function on $S^{2n-1}$ defined by $P$. Since the functional $R$ is preserved by the action of the holomorphic volume preserving transformations we deduce that $Q$ vanishes on such functions. 

At this point it is easiest to extend $Q$ in the obvious way to complex-valued functions on $S^{2n-1}$, so the conclusion above is that $Q$ vanishes on the $h^{p,1}$ and their complex conjugates $h^{1,p}$.  On each $h^{p,q}$we can compute the terms appearing in (7) explicitly.  The vector field $v$ is the generator of the circle action on the sphere so
$$ \Vert L_{v} f \Vert^{2}= (p-q)^{2}\Vert f \Vert^{2}. $$

A harmonic polynomial of degree $d$ gives an eigenfunction of the Laplacian on the sphere with eigenvalue $d(d+2n-2)$. So for $f\in h^{p,q}$ we have
$$   \Vert df\Vert^{2}= (p+q)(p+q+2n-2) \Vert f\Vert^{2}. $$

Since $df$ is the orthogonal sum of horizontal and vertical components we have
$$  \Vert d f\Vert^{2} = \Vert d_{H}f\Vert^{2} + \Vert L_{v}f\Vert^{2}$$
so $$  \Vert d_{H}f\Vert^{2}=\lambda_{pq}  \Vert f \Vert^{2}, $$ where
$\lambda_{pq}=(p+q)(p+q+2n-2) - (p-q)^{2}$.
This implies that  $\Delta_{H} f= \lambda_{pq} f$ for $f\in h^{p,q}$ so
$$   \Vert \Delta_{H} f\Vert^{2}= \lambda_{pq}^{2} \Vert f \Vert^{2}. $$
 Putting this together, if $Q$ is given by (7) we have $Q(f)= \mu_{pq} \Vert f\Vert^{2}$ for $f\in h^{p,q}$ with $p+q>0$ where 
$$  \mu_{pq}= a_{0} + a_{1}\lambda_{pq} + a_{2} (p-q)^{2} + a_{3} \lambda_{pq}^{2}. $$
One finds that the condition that $\mu_{pq}=0$ for $p=1$ or $q=1$ determines the coefficients $a_{i}$ up to an overall factor.  The solution (up to a factor) is 
$$  a_{0}= -n^{2} \ \ ,\  a_{1}=n/2\ \ \ ,\ \ \ a_{2}= (n+1)^{2}/4\ \ ,\ \ a_{3}=-1/16, $$
which gives
\begin{equation}   \mu_{pq}= -(p-1)(q-1)(p+n)(q+n). \end{equation}

 Looking a little further into the calculation one see that in fact the 
coefficient $a_{0}$ is $- 2n^{3}/((n+1)^{2} {\rm Vol}\ S^{2n-1})$, but the  crucial thing is that this is negative, so the sign of the form $Q$ on the $h^{p,q}$ is as in (8). Thus the second variation is zero on the $h^{p,1}$ and $h^{1,q}$, {\it positive} on the $h^{p,0}$ and $h^{0,q}$ and {\it negative} on the other $h^{p,q}$.

\

\
The spaces $h^{p,p}$ correspond to deformations of the sphere which are invariant under the circle action on $\bC^{n}$. For these we can obtain a simple global result.
Let $\cF_{S^{1}}$ be the space of $S^{1}$-equivariant embeddings of $S^{2n-1}$ as a pseudoconvex hypersurface in $\bC^{n}$ and let $\cF^{0}_{S^{1}}$ be the connected component containing the standard embedding.
\begin{prop}

The standard sphere is a global maximum of the functional $R$ on $\cF^{0}_{S^{1}}$.
\end{prop}

We show first that the image of any embedding  in $\cF^{0}_{S^{1}}$ is the graph of an $S^{1}$-invariant positive function  on $S^{2n-1}$, just as we considered above for small variations. In other words, we can reduce to the case of embeddings 
$m_{f}(\xi)= e^{f(\xi)} \xi$, for $\xi\in S^{2n-1}$. 

Let $m:S^{2n-1}\rightarrow \bC^{n}$ be a map in $\cF_{S^{1}}$. Clearly the image lies in $\bC^{n}\setminus\{0\}$ and we get an induced map $\underline{m}:\bC\bP^{n-1}\rightarrow \bC\bP^{n-1}$. Then the image of $m$ has the required form if and only if
 $\underline{m}$ is a diffeomorphism. In turn this is equivalent to the fact that $\underline{m}$ is an immersion. The CR structure on the image of $m$ defines an $S^{1}$-invariant contact structure $H_{m}\subset TS^{2n-1}$ and we can fix an $S^{1}$-invariant contact $1$-form $\theta_{m}$. A moment's thought shows that $\underline{m}$ is an immersion if and only if the orbits of the $S^{1}$ action are everywhere transverse to $H_{m}$. Let $\theta$ be the standard contact $1$-form on $S^{2n-1}$ with kernel $H$. We can write $\theta_{m}= g \theta+ \lambda$ where $\lambda$ vanishes on the generator of the $S^{1}$ action. So $\underline{m}$ is an immersion if and only if $g$ is nowhere zero. On the other hand there is no point where $g$ and $dg$ both vanish since at such a point
$$ \theta_{m}\wedge (d\theta_{m})^{n-1}= \lambda\wedge (d_{H}\lambda)^{n-1}=0, $$
contradicting the contact condition.

Now suppose that we have a $1$-parameter family $m_{t}$ in $\cF_{S^{1}}$, with $t\in [0,1]$, and hence a $1$-parameter family of functions $g_{t}$.
Suppose that $g_{t}$ is nowhere zero for small $t$. If $t_{0}$ were the first time when $g_{t}$ has a zero at some point $p$ then the derivative of $g_{t}$ at $p$ must also vanish and  by the discussion above this cannot happen. So we deduce that in the connected component of the standard embedding we always have $g$ nowhere zero which establishes the required statement.
(The authors do not know if there is any other component of $\cF_{S^{1}}$; this seems an interesting question.)

From here we can assume   that the map $m$ is simply given by $m(\xi) = e^{f(\xi)} \xi$ where $f$ is an $S^{1}$-invariant function on the sphere, or equivalently a function on $\bC\bP^{n-1}$. The condition that $M$ be pseudoconvex goes over to the condition that $f$ is a K\"ahler potential {\it i.e.} that $\omega_{f}>0$ for $\omega_{f}=\omega_{FS}+ i\db\partial f$, where $\omega_{FS}$ is the Fubini-Study form. 

One finds easily that 
$$A(M)=  \int_{\bC\bP^{n-1}}   \left(\omega_{f}^{n-1}\right)^{1/n+1} \left( e^{2nf}\omega_{FS}^{n-1})\right)^{n/n+1}, $$
while
$$ V= (2n)^{-1} \int_{\bC\bP^{n-1}} e^{2nf} \omega_{FS}^{n-1}. $$

Applying H\"older's inequality  with exponents
$(1/n+1), (n/n+1)$ to the first equation  we get
$$  A(M)\leq \left( \int_{\bC\bP^{n-1}} \omega_{f}^{n-1}\right)^{1/n+1}\left( \int_{\bC\bP^{n-1}} e^{2nf} \omega_{FS}^{n-1}\right)^{n/n+1}, $$
with equality if and only $\omega_{f}^{n-1}= e^{2nf}\omega_{FS}^{n-1}$. The first term on the right hand side is independent of $f$ and the second is a multiple of $V^{n/n+1}$. This implies that the functional $R$ is maximised when $f=0$. The equation characterising the maxima is the K\"ahler-Einstein equation for $\omega_{f}$ so the solutions are equivalent to the standard solution $f=0$ under projective transformations (by the Bando-Mabuchi uniqueness theorem). In other words the maximisers in ${\cal F}_{S^{1}}$ are equivalent to the sphere under linear transformations of $\bC^{n}$. 

\

{\bf Remark}

The whole discussion extends to the case of a  general Fano manifold $X$, taking $Z$ to be the complex cone defined by a power of the anticanonical bundle of $X$. While this is in general singular the functionals $A$ and $V$ are defined and the problem of maximising $R$ among suitable $S^{1}$-equivariant embeddings of a circle bundle over $X$ gives another variational formulation of the problem of finding a K\"ahler-Einstein metric on $X$.

\

The deformations defined by  the $h^{p,0}$ have a simple geometric interpretation. Let $f$ be a holomorphic function on the ball, smooth up to the boundary. Then on the boundary sphere we can write
$   f= \sum \overline{z}_{i} v_{i} $ where $v_{i}= z_{i} f$. Then $\sum v_{i} \frac{\partial}{\partial z_{i}}$ is a holomorphic vector field on the ball whose normal component at the boundary is equal to the real part of $f$. So the deformations of the sphere corresponding to the $h^{p,0}$ are effected by infinitesimal holomorphic transformations. 
Let $F:B^{2n}\rightarrow \bC^{n}$ be a holomorphic diffeomorphism to its image $\Omega \subset \bC^{n}$ which does not necessarily  preserve the volume form. Let $M=\partial \Omega= F(S^{2n-1})$. Then we can compute $R(M)$ by using the fixed domain $B^{2n}$ but the modified holomorphic $n$-form
$F^{*}(\Psi)$. Write this as $e^{G/2}\Psi$ for a holomorphic function $G$ on $B^{2n}$ and let $g$ be the real part of $G$. So $g$ is a pluriharmonic function on the ball. We find that
$$  A= \int_{S^{2n-1}}  e^{ng/n+1} , $$
while
$$   V= \int_{B^{2n}} e^{g}. $$

We thus get an interesting functional $R(g)= A/V^{n/n+1}$  on the space of pluriharmonic functions $g$ on $B^{2n}$ (smooth up to the boundary). For any such function $g$ the maximum is attained on the boundary and the outward normal derivative at a maximum point is nonnegative. We say that $g$ is \lq\lq generic'' if there is a unique maximum point and the normal derivative is strictly positive. Then straightforward estimates show that for a generic pluriharmonic function
$$   R(kg)\rightarrow\infty $$
as $k\rightarrow \infty$. It seems likely that this remains true without the generic hypothesis and it is plausible that $g=0$ is an absolute minimiser of $R$.   More generally consider the equivalence relation on pseudoconvex hypersurfaces in $\bC^{n}$ defined by  $M_{1}\sim M_{2}$ if there is a holomorphic diffeomorphism taking $M_{1}$ to $M_{2}$. Then we can seek to {\it minimize} $R$ over each equivalence class. Assuming such a minimum exists we can then seek to {\it maximize} this minimum over the set of equivalence classes. It is possible that this  minimax characterises the standard $S^{2n-1}\subset \bC^{n}$. At least this picture  is true locally, for submanifolds close to the standard sphere, by our second variation calculation.

\section{Symplectic geometry of exact forms}

There is some interesting infinite-dimensional geometry associated to the material  we have considered above and in our previous paper \cite{kn:DL}. In this section we will restrict attention to the case ${\rm dim}\ M=5$, although many of the constructions extend to general dimensions. 

Our starting point is a symplectic structure on the space $\bE$ of exact $3$-forms  on a compact oriented $5$-manifold $M$. For two such forms $\psi_{1}, \psi_{2}$ we choose $\sigma\in \Omega^{2}(M)$ such that $\psi_{1}=d\sigma$ and define
\begin{equation}  \langle \psi_{1}, \psi_{2}\rangle = \int_{M} \sigma\wedge \psi_{2} . \end{equation}
Stokes' theorem shows that this is well-defined and skew-symmetric and it endows the vector space $\bE$  with a natural symplectic structure, preserved by the action of the group ${\rm Diff}\ M$ of oriented diffeomorphisms. 
Now let $\cU\subset \bE$ be the open subset of pseudoconvex exact forms, in the sense of \cite{kn:DL}. Recall that any $\phi \in \cU$ has a canonical
factorisation
$\phi=\alpha \wedge \theta$ and Reeb vector field $v=v_{\phi}$. This defines the volume form $$\nu=\nu(\phi)= \frac{1}{2} \alpha^{2}\wedge \theta.$$
When $M$ is a pseudoconvex submanifold of a Calabi-Yau threefold and $\phi$ is the restriction of the real part of the holomorphic form this agrees with our previous definition. The integral of $\nu_{\phi}$ yields  a volume functional  $A$ on $\cU$. Using the symplectic structure this defines a Hamiltonian vector field on $\cU$ and  a flow preserving the symplectic structure.  In this setting the first variation of $A$ can be written as
$$ \delta A= \langle L_{v}\alpha \wedge \theta, \delta \phi \rangle. $$ 
In other words, the infinite-dimensional vector field on $\cU$ has value
$L_{v}\alpha\wedge \theta= d \alpha$ at $\phi\in \cU$. 

\begin{prop}

For $\phi\in \cU$ let $f_{t}:M\rightarrow M$ be the flow generated by $v=v_{\phi}$. The Hamiltonian flow generated by the functional $A$ is $f_{t}^{*}(\phi)$.
\end{prop}
This follows from the facts that $L_{v}\theta=0$ and $L_{v}\alpha\wedge \alpha=0$ which imply that the Reeb field $v_{f^{*}_{t}\phi}$ is equal to $v_{\phi}$ and then
$$  \frac{d}{dt} f_{t}^{*}\phi= L_{v}\phi= L_{v}\alpha \wedge \theta. $$

Note that while each orbit of the flow on $\cU$ is an orbit of a $1$-parameter subgroup in ${\rm Diff}M$ these $1$-parameter subgroups are not the same for different $\phi\in \cU$ so the whole flow on $\cU$ is not defined by the action of ${\rm Diff}(M)$.

\

We next look for a {\it moment map}   for the action of ${\rm Diff}(M)$ on the space $\bE$ of exact $3$-forms. Recall from multilinear algebra that we have an isomorphism
$$ \Lambda^{3} T^{*}M \cong \Lambda^{2} TM \otimes \Lambda^{5}T^{*}M . $$
Combined with the contraction  $\Lambda^{2}TM \otimes \Lambda^{3}T^{*}M\rightarrow \Lambda^{1}T^{*}M$ we  get a bilinear map
$$ \Lambda^{3}T^{*}M\otimes \Lambda^{3}T^{*}M\rightarrow \Lambda^{1}T^{*}M \otimes \Lambda^{5}T^{*}M. $$
which one checks is symmetric. Thus for a $3$-form $\phi\in \Omega^{3}(M)$ we get a quadratic expression $\phi*\phi$ in $\Gamma(T^{*}M\otimes \Lambda^{5}T^{*}M)$.
This space has a natural pairing with the space of vector fields on $M$, with values in the $5$-forms,  so for any vector field $w$ we can define
$$\mu(\phi)(w) = \frac{1}{2}\int_{M} (\phi*\phi)\ (w). $$
\begin{prop}

Restricted to the exact $3$-forms, the map $\mu$ is an equivariant moment map for the action of ${\rm Diff}\ M$, with respect to the symplectic form $\langle\ ,\ \rangle$.
\end{prop}

 The map is clearly equivariant. The moment map condition is
 that, for any vector field $w$, 
 $$  \delta \mu(w)= \langle \delta \phi, L_{w}\phi \rangle. $$
 Since $L_{w}\phi= d( i_{w}\phi)$ the right hand side is
 $$  -\int_{M} \delta\phi \wedge i_{w}(\phi). $$
 Then the statement follows from the pointwise identity
 $$  \delta(\phi*\phi)(w)= 2 \delta\phi \wedge i_{w}(\phi). $$

 Now restrict to the open subset $\cU$ of exact pseudoconvex $3$-forms.
 One finds that for $\phi=\alpha\wedge \theta$
 $$  \phi*\phi= \theta \otimes \left( \alpha^{2}\wedge \theta\right) \in T^{*}M\otimes \Lambda^{5}T^{*}M. $$
By the definition of the factorisation this can also be written as
$$  \phi *\phi = \theta \otimes \left( (d\theta)^{2}\wedge \theta\right). $$
One sees then that pseudoconvex forms $\phi_{1}, \phi_{2}$ have the same image under the moment map if and only if they define the same contact $1$-form $\theta$.  A corollary is that the Hamiltonian flow defined by any ${\rm Diff}\ (M)$-invariant functional on $\cU$ preserves the contact $1$-form, as we have seen in the particular case of the functional $A$.

Now consider the action of ${\rm Diff}\ (M)$ on the vector space $\bE\times \bE$ with the symplectic form   $\Omega_{1}$ equal to $\langle \ ,\ \rangle$ on the first factor and $-\langle\ ,\  \rangle$ on the second factor. The moment map for the action is
$$\mu_{1}(\phi_{1},\phi_{2})  = \mu(\phi_{1})-\mu(\phi_{2}).$$
Restrict to the subset $\cU\times \cU\subset \bE\times \bE$.
It follows from the remark above the zeros of $\mu_{1}$ consist of pairs
$(\alpha\wedge\theta,\beta\wedge\theta)$ with $\alpha^{2}\wedge \theta=\beta^{2}\wedge\theta=d\theta^{2}\wedge\theta$.
There is another symplectic form $\Omega_{2}$ on $\bE\times \bE$ given by
$$   \Omega_{2}\left( (\delta_{1}\phi_{1},\delta_{1}\phi_{2}), (\delta_{2}\phi_{1},\delta_{2}\phi_{2})\right)=
\langle \delta_{1}\phi_{1},\delta_{2}\phi_{2}\rangle- \langle \delta_{1}\phi_{2}, \delta_{2}\phi_{1}\rangle. $$
If we identify $\bE\times \bE$ with  the complexified vector space $\bE\otimes \bC$ ({\it i.e.} the exact complex-valued $3$-forms), then $\Omega_{1}+ i\Omega_{2}$ is the complexification of the the symplectic form $\langle\ ,\  \rangle$ on $\bE$. So $\cU\times \cU$ has the structure of an infinite-dimensional complex symplectic manifold. The moment map for the action of ${\rm Diff}\ (M)$ with respect to the real form $\Omega_{2}$ is
$$\mu_{2}(\phi_{1}, \phi_{2})(w)= 2 \int_{M} i_{w}\phi_{1}\wedge \phi_{2}. $$
The infinite-dimensional group ${\rm Diff}\ (M)$ does not have a complexification but its Lie algebra does and this complexified Lie algebra acts on $\cU\times \cU$. The combination $\mu_{\bC}= \mu_{1}+ i\mu_{2}$ is the moment map for the complexfied action with respect to the complex symplectic form. (The definition of a moment map only requires a  Lie algebra action.)
\begin{prop}
The zero set of $\mu_{\bC}$ in $\cU\times \cU$ consists of pairs with canonical factorisations $(\alpha\wedge \theta, \beta\wedge \theta)$ where
$\alpha^{2}=\beta^{2}$ and $\alpha\wedge\beta=0$.
\end{prop}

For we have  seen that the zero set  of $\mu_{1}$ consists of pairs $\alpha\wedge\theta,\beta\wedge\theta$ with $\alpha^{2}=\beta^{2}$. For such pairs the formula for $\mu_{2}$ reduces to $2\theta\otimes (\alpha\wedge\beta\wedge \theta) $,  so the vanishing of $\mu_{2}$ is the orthogonality condition $\alpha\wedge\beta=0$. 

The  orthonormal condition for $\alpha,\beta$ in Proposition 6 is exactly the condition  studied in \cite{kn:DL}. For such a pair we defined a finite-dimensional obstruction space $\cH$ and if $\cH$ vanishes then for any small deformation of $\alpha\wedge \theta$ there is a unique small deformation of $\beta\wedge \theta$ preserving the conditions. Moreover, if $(\alpha+ i\beta)\wedge \theta$ is the pull-back of the standard complex $3$-form $\Psi=dz_{1}dz_{2}dz_{3}$ on $\bC^{3}$ by  an embedding $F:M\rightarrow \bC^{3}$ (with image the boundary of a pseudoconvex domain), then the small deformations of $\alpha,\beta$ are realised by small deformation of $F$. To simplify our discussion here let us now suppose that we have an open subset $\cU_{0}\subset \cU$ on which the obstruction space $\cH$ vanishes and that for each $\alpha\wedge\theta\in \cU_{0}$ there is a unique
$\beta\wedge \theta\in \cU_{0}$ satisfying the orthonormality condition. In other words, if we write $\cL$ for the zeros of $\mu_{\bC}$ in $\cU_{0}\times \cU_{0}$ then $\cL$ is the graph of a smooth map from $\cU_{0}$ to $\cU_{0}$. On the other hand there is an open set $\cF$ of embeddings of $M$ in $\bC^{3}$ such that $\cL$ is the image of $\cF$ by an embedding $\iota$ which takes $m:M\rightarrow \bC^{3}$ to the real and imaginary parts of $m^{*}(\Psi)$.

 The general properties of moment maps mean that $\cL$ is a complex submanifold of $\cU_{0}\times \cU_{0}$, preserved by the infinitesimal action of the complexified Lie algebra. Thus for each point of $\cL$ we have a complex-linear map
$$\rho: {\rm Vect}\ (M)\otimes \bC \rightarrow T\cL. $$
 
 \begin{prop}
 \begin{itemize}
 \item $\rho$ is surjective;
 \item $\cL$ is a complex Lagrangian submanifold of $\cU_{0}\times \cU_{0}$.
 \end{itemize}
 \end{prop}
  To establish the first item we show that $\rho$ maps onto the image of the derivative $D\iota$ of $\iota$. (Recall that, under the assumptions above, this image is the tangent space $T\cL$.)
The linear map $D\iota$ takes a section $s$ of the pull-back by $m$ of the tangent bundle of $\bC^{3}$ to $d( i_{s}(\Psi)$ and this  is complex-linear. We can write the pull-back of the tangent bundle as $ TM \oplus \bR (Iv)$, so any section $s$ can be written as $w+ f Iv$ for a vector field $w$ on $M$ and real-valued function $f$. Then 
$D\iota(s)= \rho( w+ i f v)$.  

 The first item implies that $\cL$ is a complex isotropic submanifold ({\it i.e.} the restriction of the complex symplectic form to $\cL$ vanishes). This follows because the equivariant moment map
property implies that
 $$ \Omega_{\bC}( \rho(\xi_{1}),\rho(\xi_{2}))= \mu_{\bC}( [\xi_{1}, \xi_{2}])$$
 and $\mu_{\bC}$ vanishes at the given point. To show that $\cL$ is Lagrangian we have to see that the tangent space at each point is a maximal isotropic submanifold. This follows easily from the fact that $\cL$ is a graph and the two factors in $\bE\times \bE$ are Lagrangian with respect to $\Omega_{2}$.

Under our hypotheses we have a functional $V$ on $\cU_{0}$ defined by the volume of the domain with boundary $m(M)$. The first variation of $V$ is 
\begin{equation}  \delta V= -2 \langle \beta\wedge\theta, \delta\phi\rangle . \end{equation}
Thus the Hamiltonian flow on $\cU_{0}$ generated by $V$ takes
$\alpha\wedge \theta$ to $$ \cos 2t \alpha\wedge \theta- \sin 2t \beta\wedge \theta. $$
It corresponds (up to the factor $-2$) to the standard $S^{1}$ action on $\bE\otimes \bC$, restricted to $\cL$, under the projection from $\cL$ to $\cU_{0}$. We have
\begin{prop}
The functionals $A,V$ on $\cU_{0}$ Poisson commute.
\end{prop}

This can be seen from the fact that the two flows commute. Alternatively, the Poisson bracket is
$$  \{ A, V\}= \langle d\alpha, \beta\wedge\theta\rangle= \int_{M} \alpha\wedge\beta \wedge \theta=0. $$

The first variation formula (10) can be viewed in the following way. We use the symplectic form $\langle\ ,\  \rangle$ to identify $\bE$ with $\bE^{*}$ and hence $\bE\times \bE$ with $\bE\times \bE^{*}= T^{*} \bE$. Under this identification the symplectic form $\Omega_{2}$ becomes the canonical symplectic form on $T^{*}\bE$. From this point of view (10) is the statement that $\cL$ is the graph of the derivative of the function $V$ on $\cU_{0}\subset \bE$. Now use the projection from $\cL$ to $\cU_{0}$ to define a complex structure on $\cU_{0}$. The fact that $\cL$ is complex Lagrangian implies that the symplectic form $\langle\ ,\  \rangle$ is of type $(1,1)$ with respect to this complex structure--it defines an indefinite K\"ahler form. The function $V$ is a K\"ahler potential, i.e. the form is $ dId V $.

It is tempting to think of the space $\cL$ as a homogeneous space for a formal complexification of ${\rm Diff}\ (M)$ but it seems hard to give any geometric interpretation to this idea. Things work better if we consider the smaller group of contact transformations. Thus we fix a contact structure
$H\subset TM$ and write ${\rm Diff}^{H}(M)$ for the group of diffeomorphisms preserving $H$. Let $\cF_{H}$ be the set of pseudoconvex embeddings $m:M\rightarrow \bC^{3}$ which map $H$ to the complex tangent space of the image. Then ${\rm Diff}^{H}\ (M)$ acts freely on $\cF_{H}$ by composition on the right. 

\begin{prop}

$\cF_{H}$ is a complex submanifold of the space of maps from $M$ to $\bC^{3}$ and the tangent bundle  of $\cF_{H}$ is the complexification of the subbundle of tangents to the ${\rm Diff}^{H}(M)$ orbits. 
\end{prop}

Let $\theta$ be any contact $1$-form on $M$ with kernel $H$. The condition that a map $m:M\rightarrow \bC^{3}$ lies in $\cF_{H}$ can be expressed as
\begin{equation} \theta\wedge m^{*}({\rm Re}(\Psi))=0.\end{equation} At a given $m_{0}$ in $\cF_{H}$ let $\theta\wedge \alpha$ be the canonical factorisation of $m_{0}^{*}({\rm Re}(\Psi))$, and $\theta\wedge \beta$ of the imaginary part. As above, an infinitesimal deformation of $m_{0}$ has the form  $f v + g (Iv) + \xi$ where $\xi$ is section of $H$. The differentiated form of the defining condition (11) is
$$  \theta\wedge d( f \alpha + g\beta + i_{\xi} \alpha\wedge \theta) =0, $$
which is
\begin{equation}  d_{H}(f)\wedge \alpha + d_{H} g\wedge \beta + i_{\xi}(\alpha)\wedge \omega=0. \end{equation}

For any $f,g$ this equation (12) defines a unique $\xi$, so the tangent space of $\cF_{H}$ can be identified with the complex-valued functions $f+ig$ on $M$. The real valued functions, with $g=0$, give  the tangents to the orbit of
the ${\rm Diff}^{H}(M)$ action. To complete the proof of the Proposition we have to see that multiplying $(f+ig)$ by $i$ corresponds to multiplying $\xi$ by $I$--the complex structure on $H$ defined by $m_{0}$. This follows from the fact that if $ i_{\xi}(\alpha)\wedge \omega = a\wedge \alpha$ for $a\in H^{*}$ then $i_{I\xi}(\alpha)\wedge \omega = a\wedge \beta$.

The space $\cF_{H}$ is not a group but Proposition 9 implies that it carries the infinitesimal structure of one---a complexification ${\rm Diff}^{H}_{c}(M)$---in the sense discussed in \cite{kn:D1} and many other references. From this point of view the  space of unparametrised pseudoconvex hypersurfaces in $\bC^{3}$ diffeomorphic to $M$ becomes the quotient $${\rm Diff}^{H}_{c}(M)/{\rm Diff}^{H}(M). $$ One can go on to discuss (affine) geodesics in this space, which give an interpretation of $1$-parameter subgroups in ${\rm Diff}^{H}_{c}(M)$ but we will not do that here.

 Returning to our space of pseudoconvex $3$-forms on $M$, let $\cU_{0}^{H}\subset \cU_{0}$ be the space of pseudoconvex forms $\alpha\wedge \theta\in \cU_{0}$ with ${\rm ker}\ \theta=H$. Then ${\rm  Diff}^{H}(M)$ acts on $\cU_{0}^{H}$.
It also acts on the subset  $\cL^{H}= \left(\cU_{0}^{H}\times \cU_{0}^{H}\right)\cap \cL$ and, just as we saw for $\cL$, the tangent space of $\cL^{H}$ is generated by the complexified Lie algebra action. But now we have an interpretation of the \lq\lq action'' of ${\rm Diff}^{H}_{c}(M)$ on $\cL^{H}$: it is given by the restriction of the map $\iota$ from $\cF_{H}$ to $\cL^{H}$ which takes $m$ to $m^{*}(\Psi)$.


\begin{thebibliography}{99}
 \bibitem{kn:Blaschke} W. Blaschke, {\em {\"U}ber affine Geometrie I: Isoperimetrische Eigenschaften von Ellipse und Ellipsoid} Ber. Verh. Sachs. Akad. Wiss., Math.-Phys. Kl. 68 (1916), 217-239
\bibitem{kn:Deicke}  A. Deicke
{\em \"Uber die Finsler-R\"aume mit $A_i=0$} 
Arch. Math. 4, 45-51 (1953).
 \bibitem{kn:D1} S. Donaldson  {\em Symmetric spaces, K\"ahler geometry and Hamiltonian 
dynamics} In: Northern-California Seminar on 
Symplectic Geometry (Eds. Eliashberg, Weinstein)  Amer. Math. Soc.
  1999    13-34             
\bibitem{kn:DL} S. Donaldson and F. Lehmann {\em Closed 3-forms in five dimensions and embedding problems}  arXiv:2210.16208 
\bibitem{kn:Nakajima} S. Nakajima {\em \"Uber die Isoperimetrie der Ellipsoide und Eifl\"achen mit konstanter
mittlerer Affinkr\"ummung im $(n+1)$-dimensionalen Raume} 
Japanese Journ. of Math. 2, (1926) 193-196 
\bibitem{kn:Santalo} L. Santal\'o {\em An affine invariant for $n$-dimensional convex bodies} 
Port. Math. 8, (1949) 155-161 


\end{thebibliography}
\end{document}